\documentclass[12pt,twoside]{article}

\usepackage{amsmath,amssymb,amsthm,amscd,hyperref,graphicx,color,accents}
\numberwithin{equation}{section}

\newtheorem*{theorem*}{Theorem}

\def\C{\mathbb{C}}

\def\dd{\,{\rm{d}}}

\def\N{\mathbb{N}}
\def\R{\mathbb{R}}

\begin{document}

\title{\bf\large 
A Generalization of an Integral Arising in the Theory of Distance Correlation}

\author{
{Johannes Dueck,}\thanks{Institut f\"ur Angewandte Mathematik, Universit\"at Heidelberg, Im Neuenheimer Feld 294, 69120 Heidelberg, Germany.} 
\ { Dominic Edelmann,}$^*$ 
{and Donald Richards}\thanks{Department of Statistics, Pennsylvania State University, University Park, PA 16802, U.S.A.
}
}



\maketitle


\begin{abstract} 
We generalize an integral which arises in several areas in probability and statistics and which is at the core of the field of distance correlation, a concept developed by Sz\'ekely, Rizzo and Bakirov (2007) to measure dependence between random variables.

Let $m$ be a positive integer and let ${\cos_m}(v)$, $v \in \mathbb{R}$, be the truncated Maclaurin expansion of ${\cos}(v)$, where the expansion is truncated at the $m$th summand.  For $t, x \in \mathbb{R}^d$, let $\langle t,x\rangle$ and $\|x\|$ denote the standard Euclidean inner product and norm, respectively.  We establish the integral formula: For $\alpha \in \mathbb{C}$ and $x \in \mathbb{R}^d$, 
$$
\int_{{\mathbb{R}}^d}\frac{\cos_m(\langle t,x\rangle) - \cos(\langle t,x\rangle)}{\|t\|^{d+\alpha}} \,\dd t = C(d,\alpha) \, \|x\|^{\alpha},
$$
with absolute convergence if and only if $2(m-1) < \Re(\alpha) < 2m$.  Moreover, the constant $C(d,\alpha)$ does not depend on $m$.



\smallskip
\noindent
{{\em Running head}: An Integral Arising in Distance Correlation.}
\end{abstract}

\section{Introduction}
\label{introduction}
\numberwithin{equation}{section}

Sz\'ekely, Rizzo and Bakirov (2007) and Sz\'ekely and Rizzo (2009) introduced the powerful concepts of distance covariance and distance correlation as measures of dependence between collections of random variables.  In later papers, Rizzo and Sz\'ekely (2010, 2011) and Sz\'ekely and Rizzo (2012, 2013, 2014) gave applications of the distance correlation concept to several problems in mathematical statistics.  In recent years, there have appeared an enormous number of papers in which the distance correlation coefficient has been applied to many fields.  In particular, the concept of distance covariance has been extended to abstract metric spaces (Lyons, 2013) and has been related to machine learning (Sejdinovic, Sriperumbudur, Gretton, and Fukumizu, 2013); and there have been applications to the analysis of wind data (Dueck, Edelmann, Gneiting, and Richards, 2014); to detecting associations in large astrophysical databases (Mart\'inez-Gomez, Richards, and Richards, 2014) and to interpreting those associations (Richards, Richards, and Mart\'inez-Gomez, 2014); to measuring nonlinear dependence in time series data (Zhou, 2012); and to numerous other fields.  

Sz\'ekely and Rizzo (2005), in developing the foundations of distance correlation, derived an intriguing multidimensional singular integral.  It is this integral which is the subject of the present paper; and in stating the result, we denote for $t, x \in \R^d$ the standard inner product and Euclidean norm by $\langle t,x\rangle$ and $\|t\|$, respectively.  

Suppose that $\alpha \in \C$ satisfies $0 < \Re(\alpha) < 2$.  Sz\'ekely and Rizzo (2005) proved that, for all $x \in \R^d$, 
\begin{equation}
\label{sz_integral}
\int_{\R^d}\frac{1-\cos(\langle t,x\rangle)}{\|t\|^{d+\alpha}} \,\dd t 
= C(d,\alpha) \, \|x\|^{\alpha},
\end{equation}
where 
\begin{equation}
\label{normalization_constant}
C(d,\alpha)=\frac{2\pi^{d/2} \, \Gamma(1-\alpha/2)}{\alpha \, 2^{\alpha} \, \Gamma\big((d+\alpha)/2\big)}.
\end{equation}
Sz\'ekely and Rizzo defined the integral (\ref{sz_integral}) by means of a regularization procedure, where the integrals at $0$ and at $\infty$ are in a principal value sense: $\lim_{\epsilon \to 0} \int_{\R^d \setminus \{\epsilon B + \epsilon^{-1}B^c\}}$, where $B$ is the unit ball centered at the origin in $\R^d$ and $B^c$ is the complement of $B$.  

In this paper, we generalize the integral (\ref{sz_integral}) by inserting into the integrand a truncated Maclaurin expansion of the function $\cos(\langle t,x\rangle)$.  We show that the generalization is valid for all $\alpha \in \C$ such that $2(m-1) < \Re(\alpha) < 2m$, where $m$ is any positive integer.  Moreover, we prove that the generalization converges absolutely under the stated condition on $\alpha$; as a consequence, we deduce that (\ref{sz_integral}) converges without the need for regularization.  

We note that the integral (\ref{sz_integral}) arises in other areas of probability and statistics.   Indeed, in the area of generalized random fields, (\ref{sz_integral}) provides the spectral measure of a power law generalized covariance function, which corresponds to fractional Brownian motion; see Reed, Lee and Truong (1995) or Chil\`es and Delfiner (2012, p. 266, Section 4.5.6).  In mathematical analysis, a related integral is treated by Gelfand and Shilov (1964, pp. 192--195), and a similar singular integral arises in Fourier analysis in the derivation of the norms of integral operators between certain Sobolev spaces of functions (Stein, 1970, pp. 140 and 263).  

We remark that the extension of (\ref{sz_integral}) to more general values of $\alpha$ raises the intriguing possibility that a general theory of distance correlation can be developed for values of $\alpha$ outside the range $(0,2)$.

\section{The Main Result}
\label{mainresult}
\numberwithin{equation}{section}

Let $m \in \N$, the set of positive integers.  Also, for $v \in \R$, define 
\begin{equation}
\label{truncated_cosine}
\cos_m(v) := \sum_{j=0}^{m-1} (-1)^j \frac{v^{2j}}{(2j)!}
\end{equation}
to be the truncated Maclaurin expansion of the cosine function, where the expansion is halted at the $m$th summand.  

The following result generalizes (\ref{sz_integral}) to arbitrary $m \in \N$.  

\smallskip

\begin{theorem*} 
Let $m \in \N$ and $x \in \R^d$.  For $\alpha \in \C$, 
\begin{equation}
\label{main.integral}
\int_{\R^d}\frac{\cos_m(\langle t,x\rangle) - \cos(\langle t,x\rangle)}{\|t\|^{d+\alpha}} \,\dd t = C(d,\alpha) \, \|x\|^{\alpha},
\end{equation}
with absolute convergence if and only if $2(m-1) < \Re(\alpha) < 2m$, where $C(d,\alpha)$ is given in (\ref{normalization_constant}).  
\end{theorem*}

\noindent{\sc Proof.}  We shall establish the proof by induction on $m$.  

Throughout the proof, we let $B_a = \{x \in \R^d : \|x\| < a\}$ denote the ball which is centered at the origin and which is of radius $a$.  

Consider the case in which $m = 1$.  In this case, observe that for $t \in B_a$ where $a$ is sufficiently small, the function 
$$
t \mapsto \cos_1(\langle t,x\rangle) - \cos(\langle t,x\rangle) \equiv 1 - \cos(\langle t,x\rangle)
$$
is asymptotic to $\|t\|^2$.  Then the integrand in (\ref{main.integral}), when restricted to $B_a$, is asymptotic to $\|t\|^{-d-\alpha+2}$.  By a transformation to spherical coordinates to compute the integral over the unit ball $B$ we deduce that the integrand is integrable over $B_a$, and hence integrable over any compact neighborhood of the origin, if and only if $\Re(\alpha) < 2$.  

For $\|t\| \to \infty$, we apply the bound $|1- \cos(\langle t,x\rangle)| \leq 2$ to deduce that the integrand in (\ref{main.integral}) (with $m = 1$) is integrable over $\R \setminus B_a$ if and only if $\Re(\alpha) > 0$.  Consequently, for $m = 1$, the integral converges for all $x \in \R^d$ if and only if $0 < \Re(\alpha) < 2$.  

To conclude the proof for the case in which $m = 1$, we proceed precisely as did Sz\'ekely, et al. (2007, p. 2771) to obtain the right-hand side of (\ref{main.integral}).  
 
Next, we assume by inductive hypothesis that the assertion holds for a given positive integer $m$.  Note that the right-hand side of (\ref{main.integral}), as a function of $\alpha \in \C$, is meromorphic with a pole at each nonnegative integral $\alpha$.  

By (\ref{truncated_cosine}), 
$$
\cos_{m+1}(v) = \cos_m(v) + (-1)^m \frac{v^{2m}}{(2m)!}.
$$
For fixed $a > 0$, we decompose the integral (\ref{main.integral}) into a sum of three terms: 
\begin{align}
\label{decomposition}
\int_{\R^d} \frac{\cos_m(\langle t,x\rangle) - \cos(\langle t,x\rangle)}{\|t\|^{d+\alpha}} \,\dd t =  T_1 + T_2 + T_3,
\end{align}
where 
\begin{align*}
T_1 & = \int_{B_a} \frac{\cos_{m+1}(\langle t,x\rangle)-\cos(\langle t,x\rangle)}{\|t\|^{d+\alpha}} \,\dd t, \\
T_2 & = \int_{\R^d\backslash B_a} \frac{\cos_m(\langle t,x\rangle) - \cos(\langle t,x\rangle)}{\|t\|^{d+\alpha}} \,\dd t, \\
\intertext{and }
T_3 & = \frac{(-1)^{m-1}}{(2m)!} \int_{B_a} \frac{\langle t,x\rangle^{2m}}{\|t\|^{d+\alpha}} \,\dd t.
\end{align*}
We now determine the necessary and sufficient condition on the range of $\alpha$ for which the decomposition (\ref{decomposition}) entails absolute convergence of the integral.  In so doing, we examine each term individually. 

In the case of $T_1$, we apply (\ref{truncated_cosine}) to write 
\begin{equation}
\label{cosine_remainder}
\cos_{m+1}(\langle t,x\rangle)-\cos(\langle t,x\rangle) = \sum_{j=m+1}^\infty (-1)^{j+1} \frac{\langle t,x\rangle^{2j}}{(2j)!}.
\end{equation}
Proceeding formally to interchange the integral and summation, we obtain 
\begin{equation}
\label{T1series}
T_1 = \sum_{j=m+1}^\infty \frac{(-1)^{j+1}}{(2j)!} \int_{B_a} \frac{\langle t,x\rangle^{2j}}{\|t\|^{d+\alpha}} \,\dd t.
\end{equation}
To verify that this series converges absolutely, note that 
\begin{equation}
\label{T1upperbound1j}
\int_{B_a} \frac{\langle t,x\rangle^{2j}}{\|t\|^{d+\alpha}} \,\dd t 
\end{equation}
converges absolutely for all $x \in \R^d$ if and only if $\Re(\alpha) < 2j$.  Moreover, this integral clearly is a radial function of $x$, and it can be calculated exactly by a transformation to spherical coordinates.  After evaluating the integral and inserting it in the series (\ref{T1series}), we find that the series converges absolutely for all $x \in \R^d$ if and only if $\Re(\alpha) < 2(m+1)$.  

As regards the term $T_2$ we know, by inductive hypothesis, that it converges absolutely if and only if $\Re(\alpha) > 2(m-1)$.  

To analyze the term $T_3$, we note that $T_3$ is similar to (\ref{T1upperbound1j}); hence we find that $T_3$ converges absolutely if and only if $\Re(\alpha) < 2m$.  

To complete the proof, we need to evaluate $T_3$.  Let $S^{d-1}$ be the unit sphere in $\R^d$ and, for $\omega = (\omega_1,\ldots,\omega_d) \in S^{d-1}$, let ${\rm{d}}\omega$ denote the corresponding surface measure.  We define 
$$
A_{d-1} = \int_{S^{d-1}} \omega_1^{2m} \dd \omega,
$$
a constant which can be calculated exactly but whose exact value is not needed in this context.  

Similar to (\ref{T1upperbound1j}), $T_3$ is a radial function of $x$.  Thus, by a standard invariance argument and by a transformation to spherical coordinates, $t = r\omega$, where $\omega \in S^{d-1}$ and $0 \le r \le a$, we obtain 
\begin{align}
\label{T3term}
T_3 &=\frac{(-1)^{m-1}}{(2m)!} \, A_{d-1} \, \|x\|^{2m} \, \int_0^a r^{2m-1-\alpha} \dd r   \nonumber\\ 
    &= \frac{(-1)^{m-1}}{(2m)!} \, A_{d-1} \, \|x\|^{2m} \,\frac{a^{2m-\alpha}}{2m-\alpha} .
\end{align}
Moreover, the last term in (\ref{T3term}) exists for all $\alpha \in \C$ such that $\Re(\alpha) \neq 2m$ and it is a meromorphic function of $\alpha$.  

To summarize, $T_1$ converges absolutely for $\Re(\alpha) < 2(m+1)$; $T_2$ converges absolutely for $\Re(\alpha) > 2(m-1)$; and $T_3$ converges absolutely for $\Re(\alpha) < 2m$.  Therefore, the decomposition (\ref{decomposition}) is valid for $2(m-1) < \Re(\alpha) < 2m$, and it represents an analytic function which equals $C(d,\alpha) \, \|x\|^{\alpha}$ on the strip $\{\alpha \in \C: 2(m-1) < \Re(\alpha) < 2m\}$.  Hence, by analytic continuation, we obtain for $2(m-1) < \Re(\alpha) < 2(m+1)$, $\Re(\alpha) \neq 2m$,
\begin{align}
\label{decomposition2}
C(d,\alpha) \, \|x\|^{\alpha} = T_1 + T_2 + \frac{(-1)^{m-1}}{(2m)!} \, A_{d-1} \, \|x\|^{2m} \,\frac{a^{2m-\alpha}}{2m-\alpha}.
\end{align}

Now fix $2m < \Re(\alpha) < 2(m+1)$ and let $a \to \infty$ in (\ref{decomposition2}).  It is apparent that $T_2 \to 0$ and $a^{2m-\alpha} \to 0$; therefore, for $2m < \Re(\alpha) < 2(m+1)$, we obtain 
\begin{align*}
C(d,\alpha) \, \|x\|^{\alpha} = \lim_{a \to \infty} T_1 = \int_{\R^d} \frac{\cos_{m+1}(\langle t,x\rangle) - \cos(\langle t,x\rangle)}{\|t\|^{d+\alpha}} \,\dd t,
\end{align*}
which concludes the proof.  
$\qed$

\bigskip
\medskip

In conclusion, we are intrigued by the possibility of applying (\ref{main.integral}) to develop a general theory of distance correlation for values of $\Re(\alpha) > 2$.  We expect, {\it inter alia}, that such a theory will lead for sufficiently large $\Re(\alpha)$ to distance correlation analyses of data modeled by random vectors which do not have finite first moments, e.g., the multivariate stable distributions of index less than $2$.  Moreover, although the integral  (\ref{main.integral}) diverges for $\Re(\alpha) = 2m$, our results raise the possibility of developing a theory of distance correlation at the poles by modifying (\ref{main.integral}) to attain convergence as $\Re(\alpha)$ converges to the poles.  

Finally, we remark that our decomposition (\ref{decomposition}) was motivated by the ideas of Gelfand and Shilov (1964, p. 10).  

\bigskip
\bigskip

{\bf Acknowledgments.}  The research of Dueck and Edelmann was supported by the {\it Deutsche Forschungsgemeinschaft} (German Research Foundation) within the program ``Spatio/Temporal Graphical Models and Applications in Image Analysis,'' grant GRK 1653.  The research of Richards was supported by the U.S. National Science Foundation, grant DMS-1309808; and by a Romberg Guest Professorship at the Heidelberg University Graduate School for Mathematical and Computational Methods in the Sciences, funded by the German Universities Excellence Initiative grant GSC 220/2.

\vspace{7mm}

\section*{References}

\newenvironment{reflist}{\begin{list}{}{\itemsep 0mm \parsep 0.2mm
\listparindent -5mm \leftmargin 5mm} \item \ }{\end{list}}

\vspace{-7mm}

\begin{reflist}

\parskip=2.5pt

Chil\`es, J.~P.,~and~Delfiner P.~(2012).~{\sl Geostatistics: Modeling Spatial Uncertainty}, second edition.  Wiley, New York.

Dueck, J., Edelmann, D., Gneiting, T., and Richards, D. (2014).   The affinely invariant distance correlation.  {\em Bernoulli}, {\bf 20}, 2305--2330.

Gelfand, I. M.~and Shilov, G. E. (1964).  {\sl Generalized Functions}, Volume 1.  Academic Press, New York.

Lyons, R.~(2013).  Distance covariance in metric spaces. {\em Ann. Probab.} {\bf 41}, 3284--3305.

Mart\'inez-G\'omez,~E., Richards,~M.~T., and Richards, D.~St.~P.~(2014).  Distance correlation methods for discovering associations in large astrophysical databases.  {\em Astrophys. J.}, {\bf 781}, 39 (11 pp.).

Reed, I. S., Lee, P. C., and Truong, T. K. (1995).  Spectral representation of fractional Brownian motion in $n$ dimensions and its properties.  {\em IEEE Trans. Inform. Theory}, {\bf 41}, 1439--1451.

Richards,~M.~T., Richards, D.~St.~P.,~and Mart\'inez-G\'omez, E. (2014).  Interpreting the distance correlation results for the COMBO-17 survey.  {\em Astrophys.~J.~Lett.}, {\bf 784}, L34 (5 pp.).

Rizzo, M.~L.~and Sz\'ekely, G.~J.~(2010).  DISCO analysis: A nonparametric extension of analysis of variance.  {\em Ann. Appl. Statist.}, {\bf 4}, 1034--1055.

Rizzo, M.~L.~and Sz\'ekely, G.~J.~(2011).  E-statistics (energy statistics).  {\tt R} package, Version 1.4-0, \url{http://cran.us.r-project.org/web/packages/energy/index.html}.

Sejdinovic, D., Sriperumbudur, B., Gretton, A., and Fukumizu, K. (2013).  Equivalence of distance-based and RKHS-based statistics in hypothesis testing.  {\em Ann. Statist.}, {\bf 41}, 2263--2291.

Stein, E. M. (1970).  {\sl Singular Integrals and Differentiability Properties of Functions}.  Princeton University Press, Princeton, NJ.

Sz\'ekely, G. J.~and Rizzo, M.~L.~(2005).  Hierarchical clustering {\it via} joint between-within distances: Extending Ward's minimum variance method.  {\em J. Classification}, {\bf 22}, 151--183.

Sz\'ekely, G. J.~and Rizzo, M.~L.~(2009).  Brownian distance covariance (with discussion).  {\em Ann. Appl. Statist.}, {\bf 3}, 1236--1265.

Sz\'ekely, G. J.~and Rizzo, M.~L.~(2012).  On the uniqueness of distance correlation.  {\em Statist. Probab. Lett.}, {\bf 82}, 2278--2282.

Sz\'ekely, G. J.~and Rizzo, M.~L.~(2013).  The distance correlation $t$-test of independence in high dimension.  {\em J. Multivariate Anal.}, {\bf 117}, 193--213.

Sz\'ekely, G. J.~and Rizzo, M.~L.~(2014).  Partial distance correlation with methods for dissimilarities.  {\em Ann. Statist.}, {\bf 42}, 2382--2412.

Sz\'ekely, G.~J., Rizzo, M.~L.,~and Bakirov, N.~K.~(2007).  Measuring and testing independence by correlation of distances.  {\em Ann. Statist.}, {\bf 35}, 2769--2794.

Zhou, Z. (2012).  Measuring nonlinear dependence in time-series, a distance correlation
approach. {\em J. Time Series Anal.}, {\bf 33}, 438--457.

\end{reflist}

\end{document}